\documentclass[12pt]{amsart}
\usepackage{amssymb}
\usepackage{latexsym}
\usepackage{psfig}
\newcommand{\R}{{\mathbb R}}
\newcommand{\C}{{\mathbb C}}
\newcommand{\HH}{{\mathbb H}}

\newtheorem{lemma}{Lemma}[section]

\newtheorem{theorem}[lemma]{Theorem}

\newtheorem{problem}[lemma]{Problem}

\begin{document}
\title{ Locating Anisotropies in
Electrical Impedance Tomography}
\author{Erkki Somersalo}
\address{\hskip-\parindent
Erkki Somersalo\\
Helsinki University of Technology\\ 
Institute of
Mathematics\\
P.O. Box 1100, FIN--02015 HUT, Finland}
\email{erkki.somersalo@hut.fi}

\thanks{This article was written during the author's visit
at Mathematical Sciences Research Institute (MSRI), Berkeley in August 
and September, 2001. Research at MSRI is supported in part by NSF grant 
DMS-9701755. The author
wishes to thank Dr. Fioralba Cakoni, prof. David Colton and prof. Rainer
Kress for useful discussions as well as MSRI for the hospitality and for 
bringing us together.}

\begin{abstract}
In this article, we consider the problem
of finding the support of an inhomogenous possibly
anisotropic inclusion in a background of constant electric
conductivity from the electrical impedance tomography
data at the boundary of a bounded body. The article discusses
the linear sampling method applied to this problem. A
practical algorithm for solving this problem is suggested.
\end{abstract}

\maketitle

\section{Introduction}

In this article, the following 
electrical impedance tomography problem (EIT) is
considered: On the surface of a body with unknown
impedance distribution, one applies
a set of prescribed electric currents
and measures the corresponding voltages on the surface.
From this information, one seeks to estimate the internal
structure of the body. Potential application areas of the EIT
range from medical imaging and monitoring to industrial
process monitoring and nondestructive material testing.

Often, the materials encountered in applications are
anisotropic, i.e., the electromagnetic properties of the
medium depend on the direction. It is well known that in
general, the anisotropic EIT problem allows no unique
solution. Hence, in the presence of anisotropies one has
either to use {\em a priori} information complementary
to that obtained from the measurements or one has to
confine to more modest goals than recovering the full
information of the material parameters. In this article,
the latter approach is taken. More precisely,
the following inverse problem is studied: Assume that
in an isotropic body with otherwise
known electric properties, there
is an unknown possibly anisotropic inclusion. Given the
EIT data on the surface of the body, estimate the {\em support}
of the unknown inclusion.

In recent years, a number of articles have been
published, where the goal is to determine the shape
of an inclusion based on either far-field or near field
measurements with various probing modalities, see e.g.
\cite{bruhl00}, \cite{cakoni01}, \cite{colton97},
\cite{colton97a}, \cite{colton00a},
\cite{ikehata99}, \cite{ikehata00}, \cite{kirsch99},
\cite{colton00}
and the references in these articles. 
The starting
point of the present work is the linear sampling
method originally introduced in the article 
\cite{colton97}. The ideas here come close
to those presented in \cite{bruhl00} and 
\cite{cakoni01}.

\section{The inverse problem}

Let $B\subset\R^n$, $n=2,3$, denote 
a bounded simply connected domain with a $C^2$ smooth connected boundary.
When low-frequecy time-harmonic electromagnetic
field is induced in the body, within the
quasi-static approximation of Maxwell's
equations the electric field is written in terms
of the voltage potential $u$ that satisfies
\[
 \nabla\cdot\gamma\nabla u = 0\mbox{ in $B$}.
\]
Here, $\gamma=\gamma(x)\in\C^{n\times n}$ is the admittance distribution
that in terms of the conductivity $\sigma
(x)\in\R^{n\times n}$,
permittivity 
$\varepsilon(x)\in\R^{n\times n}$ and frequency $\omega>0$
is given as
\[
 \gamma = \sigma -{\rm i}\omega\varepsilon.
\] 
Assume that $D\subset B$ is an open set
that has a smooth boundary and consists of one or several simply connected
components. Furthermore, assume that 
$\partial D\cap\partial B=\emptyset$. We assume that the
permittivity of the material outside $D$ is so low
that it can be neglected, i.e.,
$\gamma\approx\sigma$ in $B\setminus\overline D$.
Furthermore, we assume that the material is isotropic
outisde $D$. To simplify the discussion, we
assume that in fact, $\sigma=1$ outside $D$. It
is not hard to see that the results in this article
can be generalized in a straightforward manner to cover
the case where $\sigma$ is a strictly positive
scalar function outside $D$.
These assumptions lead to the following model:
By denoting by $\chi_D$ the characteristic 
function of the set
$D$, the admittance $\gamma$ is of the form
\[
\gamma(x) =1 +h(x)\chi_D(x).
\]
Furthermore, we assume that the perturbation
$h\in C^1(B,\C^{n\times n})$ is symmetric and satisfies in addition that
for some positive constants $\alpha$, $\beta>0$ and for all
$\zeta\in\C^n$,
\begin{eqnarray}
 {\rm Re}\big(z\overline\zeta\cdot\gamma(x)\zeta\big) 
 &\geq& \alpha |\zeta|^2,\mbox{ for all $x\in B$, for some $z\in\C$}
\label{assumption1}\\
\noalign{\vskip4pt}
 {\rm Im}\big(\overline\zeta\cdot h(x)\zeta\big) &\leq& -\beta |\zeta|^2,
\mbox{ $x$ in an open set ${\mathcal O}\subset D$}\label{assumption2}
\end{eqnarray}

To describe the measurement, we fix the following notation. By
$H^s_0(\partial B)$ (with the notation
$H^0_0(\partial B)=L^2_0(\partial B)$) 
we denote the class of Sobolev functions of smoothness
index $s$ over the boundary with the restriction
\[
 \tau f =\langle f,1\rangle = 0,
\]
the brackets denoting the natural pairing between $H^s(\partial B)$
and $H^{-s}(\partial B)$. since $\tau:H^s(\partial D)\to\C$
is continuous, $H^s_0(\partial B)$ is a closed subspace.

Assume that one applies an electric current $f\in H^{-1/2}_0(\partial B)$
on the surface of the body $B$. Then the voltage potential $u$ satisfies
\begin{eqnarray}
\nabla\cdot\gamma\nabla u &=& 0\label{forward}\\
\noalign{\vskip4pt}
\partial u\big|_{\partial B} &=& f\nonumber,
\end{eqnarray}
where we have denoted  the normal derivative of $u$ at the 
boundary as $\partial u|_{\partial B}$.
One can see, by a standard argument using the Lax-Milgram lemma and
Fredholm theory that the problem (\ref{forward}) has a unique solution
$u\in H^1_0(B)$, where we have denoted
\[
 H^1_0(B)=\left\{ u\in H^1(B)\mid u\big|_{\partial B}\in H^{1/2}_0
 (\partial B)\right\}.
\]
For later reference, let us denote by $T$ the continuous solution
operator
\[
 T:H_0^{-1/2}(\partial B)\to H_0^1(B),\quad f\mapsto u.
\]
Further, we define the Neumann-to-Dirichlet map $\Lambda$ as
\[
\Lambda:H_0^{-1/2}(\partial B) \to H_0^{1/2}(\partial B),\quad
f\mapsto u|_{\partial B},
\]
where $u$ is the unique solution of the problem
(\ref{forward}). We shall also consider the forward
problem when no inclusion is present, i.e., the
boundary value problem
\begin{eqnarray}
\Delta v &=& 0\label{forward0}\\
\noalign{\vskip4pt}
\partial v\big|_{\partial B} &=& f\nonumber,
\end{eqnarray}
and the corresponding solution operator and Neumann-to-Dirichlet maps,
\[
 T_0:H_0^{-1/2}(\partial B)\to H_0^1(B),\quad f\mapsto v,
\] 
\[
\Lambda_0:H_0^{-1/2}(\partial B) \to H_0^{1/2}(\partial B),\quad
f\mapsto v|_{\partial B}.
\]
The inverse problem studied in this article can be formulated as follows:

\begin{problem}
Given  the Neumann-to-Dirichlet 
map $\Lambda$,
determine the support $D$ of the perturbation. 
\end{problem}

In the following section, we discuss the linear
sampling method that gives a practical way of
estimating the support.

\section{The linear sampling method}

The linear sampling method discussed in this article is based
on the use of certain singular solutions. Therefore, let us
define the singular solution $B\ni x\mapsto\Phi(x,y,\hat\alpha)$, 
where $y\in B$ is a parameter and $\hat\alpha\in\R^n$ is
a unit vector, as a solution of the following homogenous
Neumann problem,
\begin{eqnarray}
\Delta \Phi &=& 
\hat\alpha\cdot\nabla\delta(x-y)\label{singular}\\
\noalign{\vskip4pt}
\partial \Phi \big|_{\partial B} &=& 0.\nonumber
\end{eqnarray}
Physically, the singular solution corresponds to
the electromagnetic potential created by a dipole source
at $y$ pointing in the direction $\hat\alpha$. 

In terms of the operator $T_0$, we may write the singular solution
as
\[
 \Phi(x,y,\hat\alpha) =\hat\alpha\cdot\vec\Psi(x-y)
 -T_0(\vec \alpha\cdot\partial\vec\Psi(x-y))-c_y,
\]
where
\[
\vec\Psi(x-y) =-\frac 1{2^{n-1}\pi}\frac{x-y}{|x-y|^n},\quad
 n=2,3,
\]
and $c_y$ is a constant equal to the integral of $x\mapsto\alpha\cdot
\vec\Psi(x-y)$ over the boundary $\partial B$. 
We shall
use the notation $\Phi_y(x)=\Phi(x,y,\hat\alpha)$, suppressing
the dependence on the direction as this plays little role in the
discussion to ensue.

The starting point of the linear sampling method lies in the
following observation that we formulate as a lemma for later
reference.
In the following, we shall use the notation
$\phi_y(x) =\Phi_y(x)\big|_{\partial B}$. 

\begin{lemma}\label{y not in D}
Assume that $y\in B\setminus\overline D$. 
Then $\phi_y\notin{\rm Ran}(\Lambda-\Lambda_0)$.
\end{lemma}

{\em Proof.} Assume on the contrary
that $\phi_y\in{\rm Ran}(\Lambda-\Lambda_0)$, and
let $\phi_y =(\Lambda-\Lambda_0)\psi_y$, $\psi_y\in H^{-1/2}_0(\partial B)$.
Let $u_y=T\psi_y$ and $v_y=T_0\psi_y$ denote the solutions of the problems
(\ref{forward}) and (\ref{forward0}) with the boundary data $\psi_y$,
respectively. Set $w_y=u_y-v_y$. Now we observe that
\[
 \partial w_y = 0,\quad w_y\big|_{\partial B} =u_y-v_y =
 (\Lambda-\Lambda_0)\psi_y = \phi_y,
\]
i.e., the Cauchy data of $w_y$ and $\Phi_y$ coincide on $\partial B$.
By Holmgren's Uniqueness Theorem, we have
$w_y =\Phi_y$ in $B\setminus(\overline D\cup\{y\})$. The claim of the
lemma follows now, since $w_y$ has no singularity at $y$ while $\Phi_y$
is singular.\hfill$\Box$

Unfortunately, the converse is not in general true: When $y\in D$,
there is no guarantee that $\phi_y\in{\rm Ran}(\Lambda-\Lambda_0)$.
As in the case of inverse scattering problems, we have to confine
to an approximate solution of the equation $(\Lambda-\Lambda_0)\psi_y
=\phi_y$. To this end, we need to introduce some notations and
definitions, and prove a number of auxiliary results.

To get a handle of the following definition, assume for a while
that $\phi_y =(\Lambda-\Lambda_0)\psi_y$ for some $\psi_y\in
H^{-1/2}_0(\partial B)$ and $y\in D$. By denoting again
$u_y=T\psi_y$
and $v_y=T_0\psi_y$, we observe as in the proof of the previous
lemma that in $B\setminus\overline D$, $u_y-v_y =\Phi_y$,
and in $D$, 
\[
 \nabla\cdot\gamma\nabla u_y=\Delta v_y=0.
\]
At the boundary $\partial D$, the solutions must satisfy
\begin{eqnarray*}
 u_y\big|_{\partial D}^- -v_y\big|_{\partial D}^- 
&=& \Phi_y\big|_{\partial D},\\
\noalign{\vskip4pt}
\partial_\gamma u_y\big|_{\partial D}^- -\partial
v \big|_
{\partial D}^- &=& \partial\Phi_y\big|_{\partial D}.
\end{eqnarray*}
Above, the notation $\partial_\gamma u_y|_{\partial D}^- 
=n\cdot\gamma\nabla u_y
|_{\partial D}^-$ for the conormal derivative was used,
and the subscript "-" indicates that the traces are
from inside of $\partial D$.

Hence, in order to investigate the equation
$(\Lambda-\Lambda_0)\psi_y=\phi_y$, 
it is natural to study the following interior transmission
problem of impedance tomography. 

\begin{problem}\label{itp}
The interior transmission problem (ITP) of electrical impedance
tomography with boundary data $(f,g)\in H^{1/2}(\partial D)\times
H^{-1/2}(\partial D)$ is to
find functions $(u,w)\in (H^1(D)\times H^1(D))/\C$ 
satisfying
the equations
\[
 \nabla \cdot\gamma\nabla u =\Delta v =0\mbox{ in $D$},
\]
with
\begin{eqnarray*}
 u\big|_{\partial D} -v\big|_{\partial D} &=& f,\\
\noalign{\vskip4pt}
\partial_\gamma u\big|_{\partial D} -\partial\big|_
{\partial D} v &=& g.
\end{eqnarray*}
\end{problem}

Observe that above, the solution of the interior transmission problem
can be unique only up to an additive constant. The space
$(H^1(D)\times H^1(D))/\C$ consists of equivalence classes of the
relation
\[
 (u,v)\sim (u',v') \mbox{ if and only if }  (u,v) =(u'-c,v'-c),\quad c\in\C.
\]

It turns out that under the assumptions made about the admittance
$\gamma$, the ITP has a unique solution. We formulate this as a lemma.

\begin{lemma}\label{uniqueness}
Assume that the admissivity $\gamma$ satisfies the 
conditions (\ref{assumption1}) and
(\ref{assumption2}). Then 
the interior transmission problem \ref{itp} has a unique solution.
\end{lemma}

{\em Proof:} Tho show that there are at most one solution,
assume that $(u,v)$ satisfy the homogenous interior transmission
problem, i.e., $(f,g)=(0,0)$.
By integrating by parts we obtain
\begin{eqnarray*}
 0 &=& \int_{\partial D}({\overline v}\partial v 
-v\overline{\partial v})dS
= \int_{\partial D}({\overline u}\partial_\gamma u 
-u\overline{\partial_\gamma u})dS\\
\noalign{\vskip4pt}
&=& 2{\rm i}\int_D{\rm Im} 
(\overline{\nabla u}\cdot h\nabla u) dx.
\end{eqnarray*}
By the assumption (\ref{assumption2}), we obtain that $u$ is constant in
an open set of $D$, and by the unique continuation property it is a
constant in the whole of $D$. Thus $v$ has the same Cauchy data with
a constant solution on $\partial D$ and so also $v$ is constant in $D$,
which completes the proof of uniqueness.

To prove the existence, we refer to the article \cite{cakoni01}, where
the interior transmission problem for the scattering case was studied.
It is not hard to see, that the same argument goes through here as well.
We leave the details out.
\hfill$\Box$

Although the ITP is closely related to the existence of the solution of the
equation $(\Lambda-\Lambda_0)\psi_y=\phi_y$, there is no
equivalence, since we cannot in general extend the solutions $(u,v)$
of the ITP from $D$ to the whole domain $B$. Therefore, we consider
only such solutions that have the extension property. This is the motivation
for the following considerations.

Let $\Omega\subset\R^n$ denote a ball that contains $\overline B$
in its interior. Further, let $G(x,y)$ denote Green's function
of the Laplacian,
\[
 G(x,y)=\left\{\begin{array}{ll}{\displaystyle \frac 1{2\pi}{\rm log}
\frac 1{|x-y|}},& n=2,\\ & \\
{\displaystyle \frac 1{4\pi}\frac 1{|x-y|}},& n=3.
\end{array}\right.
\]
We define the potential operator
\[
 {\mathcal S}_D:L^2(\partial\Omega)\to C^\infty(D),\quad
 \omega\mapsto \int_{\partial\Omega}G(x,z)\omega(z)dS(z),\quad
x\in D.
\]
Further, we define a particular class of
harmonic functions in $D$ as
\[
\HH(D)=\{u\mid u={\mathcal S}_D\omega,\;
\omega\in L^2(\partial\Omega)\}.
\]
We also define the operators
\[
 K:L^2(\partial\Omega)\to H^{1/2}(\partial D),\quad
 \omega\mapsto\int_{\partial\Omega}G(x,z)\omega(z)dS(z),
\quad x\in \partial D,
\]
and, finally
\[
 L:L^2(\partial\Omega)\to H^{-1/2}_0(\partial B),\quad
 \omega\mapsto\partial\int_{\partial\Omega}G(x,z)\omega(z)dS(z),\quad x\in\partial B.
\]
Observe that by Green's formula, the integral of $L
\omega$ over $\partial B$ automatically vanishes.

The class $\HH(D)$ has the following approximation property of
harmonic functions in $D$.

\begin{theorem}\label{density}
For each $\varepsilon>0$ and 
$v\in H^1(D)$ satisfying $\Delta v=0$ in the weak
sense there is $v^\varepsilon\in\HH(D)$ such that
$\|v-v^\varepsilon\|_{H^1(D)}<\varepsilon$.
\end{theorem}

The proof of this theorem is based on the following density
result.

\begin{lemma}
The operator $K$ has dense range in both $L^2(\partial D)$
as in $H^1(\partial D)$.
\end{lemma}

{\em Proof:} The proof is quite similar to the corresponding one in
the article \cite{colton01}. To prove the denseness in 
$L^2(\partial D)$,
assume that $\eta\in L^2(\partial D)$ is such that for all $\omega
\in L^2(\partial\Omega)$,
\[
 (K\omega,\eta)_{L^2(\partial D)} =(\omega,K^*\eta)_{L^2(\partial
\Omega)} = 0,
\]
implying that
\[
 K^*\eta(x)=\int_{\partial D} G(y,x)\eta(y)dS(y) = 
\int_{\partial D} G(x,y)\eta(y)dS(y) = 0
\]
for $x\in\partial \Omega$. Then, the function $w$ defined as
\[
w(x)=\int_{\partial D} G(x,y)\eta(y)dS(y),\quad x\in \R^n,
\]
is harmonic both in $\R^n\setminus\overline D$ and $D$.
It has vanishing Dirichlet boundary data on $\partial \Omega$ and $w\rightarrow 0$ as $|x|\rightarrow
\infty$,
so $w=0$ in $\R^n\setminus\overline\Omega$, and by the unique continuation
principle, $w=0$ in $\R^n\setminus\overline D$. By the continuity of the
single layer potential, this implies also that $w\big|_{\partial D}^-=0$
and hence $w=0$ in $D$. The conclusion $\eta=0$ follows from
the well known jump relation $\eta = \partial w\big|_{\partial D}^-
-  \partial w\big|_{\partial D}^+$ of the normal derivatives.

To prove the denseness in $H^1(\partial D)$, we equip it with the inner
product
\[
 (\eta,\mu)_{H^1(\partial D)} =\int_{\partial D}\bigg(
 \overline{\eta(x)}\mu(x) + \overline{{\rm Grad}\eta(x)}
 \cdot{\rm Grad}\mu(x)\bigg)dS(x),
\]
where ${\rm Grad}$ is the surface gradient on $\partial D$.
By denoting by $K^\dagger$ the adjoint of $K$ as a mapping from
$L^2(\partial\Omega)$ to $H^1(\partial D)$, assume that 
we have
\[ 
 (K\omega,\eta)_{H^1(\partial D)} =(\omega,K^\dagger \eta)
_{L^2(\partial\Omega)} = 0
\]
for all $\omega\in L^2(\partial\Omega)$
or explicitly,
\[
 K^\dagger\eta(x) =\int_{\partial D}\bigg(G(y,x)\eta(y) +{\rm Grad}_y
 G(y,x)\cdot{\rm Grad}\eta(y)\bigg)dS(y)=0,
\]
when $x\in\partial \Omega$. 
To rewrite the second term in a more convenient form,
assume, for a while that $\eta$ is smooth. Then, by
Gauss' surface divergence theorem and by using
$G(x,y)=G(y,x)$, we obtain for $x\notin\partial D$,
\begin{eqnarray*}
 \int_{\partial D}{\rm Grad}_yG(y,x)\cdot
{\rm Grad}\eta(y) dS(y) &=&
 -\int_{\partial D}G(y,x)\cdot
{\rm Div}{\rm Grad}\eta(y) dS(y)\\
\noalign{\vskip4pt}
&=&\int_{\partial D}{\rm Grad}_yG(x,y)\cdot
{\rm Grad}\eta(y) dS(y),
\end{eqnarray*}
and by extension, this holds for all $\eta\in H^1
(\partial D)$. Furthermore, since
$\nabla_xG(x,y)=-\nabla_yG(x,y)$, we obtain
\begin{eqnarray*}
\int_{\partial D}{\rm Grad}_yG(x,y)\cdot
{\rm Grad}\eta(y) dS(y) &=&
\int_{\partial D}\nabla_yG(x,y)\cdot
{\rm Grad}\eta(y) dS(y)\\
\noalign{\vskip4pt}
&=&-\nabla\int_{\partial D}G(x,y)\cdot
{\rm Grad}\eta(y) dS(y).
\end{eqnarray*}
Hence, we define 
\[
w(x)=\int_{\partial D}G(x,y)\eta(y)dS(y) -\nabla\cdot\int_{\partial D}
 G(x,y){\rm Grad}\eta(y)dS(y),\quad x\in\R^n\setminus
\partial D,
\]
which is harmonic both inside and outside of $D$,
and by the above considerations, 
$w\big|_{\partial\Omega} = K^\dagger\eta=0$. As 
above, we conclude that $w=0$ in $\R^n\setminus
\overline D$.
On $\partial D$, the jump relations for vector
potentials (see \cite{colton98}, Theorem 6.12),
\[
 w\big|^\pm_{\partial D}(x) =\int_{\partial D}\bigg(G(x,y)\eta(y)-
 \nabla_xG(x,y)\cdot {\rm Grad}\eta(y))\bigg)dS(y).
\]
This expression is the adjoint of the single 
layer operator
\[
 S:L^2(\partial D)\to H^1(\partial D),\quad
\eta\mapsto\int_{\partial D}G(x,y)\eta(y)dS(y),
\]
see e.g. \cite{colton98}, (pp. 43--44) or \cite{colton01}. 
Hence, for all $\psi\in L^2(\partial D)$, we have
\[
 (S\psi,\eta)_{H^1(\partial D)}=(\psi,w\big|_{\partial D})_{L^2(\partial
 D)}=0.
\]
By the uniqueness of the interior Dirichlet problem for the
Laplacian, we now deduce that $S$ is injective, so choosing
$\psi =S^{-1}\eta$ we obtain that $\eta=0$. 
\hfill$\Box$

With the aid of the above lemma, we prove Theorem \ref{density}.

{\em Proof of Theorem \ref{density}:} Clearly, since $H^1(\partial D)$
is dense in $H^{1/2}(\partial D)$, the range of $K$ 
is also dense in
$H^{1/2}(\partial D)$. Let $v\in H^1(D)$ be harmonic and $\varepsilon
>0$ be given. We choose first $\omega\in L^2(\partial\Omega)$ such
that
\[
 \|K\omega-v\big|_{\partial D}\|_{H^{1/2}(\partial D)}<\delta
\]
for some $\delta=\delta(\varepsilon)>0$ to be defined later. Let
\[
 v^\varepsilon ={\mathcal S}_D\omega
 \in \HH(D).
\]
Then, by the continuity of the Dirichlet problem with respect
to the boundary data, we have
\[
 \|v-v^\varepsilon\|_{H^1(D)}\leq C
 \|K\omega-v\big|_{\partial D}\|_{H^{1/2}(\partial D)}<C\delta,
\]
and the claim follows by choosing $\delta =\varepsilon/C$.
\hfill$\Box$

It is clear that if we define the space $\HH(B)$ as
functions of the form $u={\mathcal S}_B\omega$ using obvious
notations, the harmonic functions in $B$ can be approximated
by $\HH(B)$--functions. This observation gives us the
following density result.

\begin{lemma}\label{density of L}
The range of the operator $L$ is dense in $H^{-1/2}_0(\partial B)$.
\end{lemma}

{\em Proof:} Let $\varepsilon>0$ and $\psi\in H^{-1/2}_0(\partial B)$
be given. Define a harmonic function $v$ in $B$ as $v=T_0\psi$.
By the observation above, we can find $v^\varepsilon
={\mathcal S}_B\omega\in\HH(B)$
with 
\[
 \|v-v^\varepsilon\|_{H^1(B)}<\delta,
\]
where $\delta=\delta(\varepsilon)$ is fixed later. But
by the trace theorem,
\[
 \|\partial v-\partial v^\varepsilon\|_{H^{-1/2}(\partial B)}
 =\|\psi -L\omega\|_{H^{-1/2}(\partial B)}
\leq C\|v-v^\varepsilon\|_{H^1(B)}<C\delta,
\]
so by choosing $\delta=\varepsilon/C$ the desired approximation
follows.\hfill$\Box$

The  counterpart of Lemma \ref{y not in D} that we want to prove
for $y\in D$ is the following.

\begin{lemma}
Assume that $y\in D$. Then the interior transmission problem
with the boundary data $(f,g)=(\Phi_y\big|_{\partial D},
\partial\Phi_y\big|_{\partial D})$
has a unique solution $(u,v)$ with $v\in\HH(D)$ if and only if
$\phi_y=(\Lambda-\Lambda_0)L\omega$ for some $\omega\in L^2
(\partial\Omega)$.
\end{lemma}

{\em Proof:} Assume that $(u,v)$ is the unique solution of the
ITP with $v={\mathcal S}_D\omega$. 
We can extend $v$ to the whole $B$
by setting $v={\mathcal S}_B\omega$, and extend $u$ to whole $B$ by
defining $u=v+\Phi_y$ in $B\setminus\overline D$. From the
ITP boundary conditions, it follows now that $u$ thus defined 
satisfies the equation $\nabla\cdot\gamma\nabla u=0$ in the
weak sense in $B$, and at the boundary $\partial B$,
\[
 \partial u\big|_{\partial B} =\partial(v+\Phi_y)\big|_{\partial B}
=\partial v\big|_{\partial B}=L \omega.
\]
Hence, we have
\[
 \psi_y=(u-v)\big|_{\partial B} =(\Lambda-\Lambda_0)L\omega.
\]

To prove the converse, let $\psi_y=(\Lambda-\Lambda_0)
L\omega$ for some $\omega\in L^2(\partial\Omega)$.
We define $u=TL\omega$ and
$v=T_0L\omega$. As in the proof of Lemma
\ref{y not in D}, we see that $u-v=\Phi_y$ in
$B\setminus\overline D$, and hence $(u,v)$ satisfy the ITP
with the boundary data $(f,g)=(\Phi_y\big|_{\partial D},
\partial\Phi_y\big|_{\partial D})$.
On the other hand, the function 
\[
v_0=v-{\mathcal S}_B\omega 
\] 
is harmonic in $B$ and at  the boundary, 
\[
\partial v_0 = \partial v -L\omega = 0,
\]
implying that $v_0$=constant. By the definition of $T_0$,
we also see that the integral of $v$ and thus $v_0$ over
the boundary vanishes, so $v_0=0$ and the claim follows.
\hfill$\Box$

The above lemma does not help us much since in general, the
unique solution of the ITP is not such that $v\in\HH(D)$. However,
we can always find an approximate solution, as
the following theorem states.

\begin{theorem}\label{y in D theorem}
Assume that $y\in D$. Then for any $\varepsilon>0$
the  equation $(\Lambda-\Lambda_0)L\omega
=\phi_y$ has an approximate solution in $\omega_y^\varepsilon\in
L^2(\partial \Omega)$, i.e., $\omega^\varepsilon_y$ is the
satisfies the estimate
\begin{equation}\label{quasisol}
\|(\Lambda-\Lambda_0) L \omega_y^\varepsilon- 
\phi_y\|_{H^{1/2}(\partial B)} <\varepsilon.
\end{equation}
Furthermore, when $y$ approaches the boundary $\partial D$,
 $\|\omega_y^\varepsilon\|_{L^2(\partial\Omega)}
\rightarrow\infty$.
\end{theorem}

{\em Proof:} Let $(u_y,v_y)$ be the solution of the interior transmission
problem with the transmission data $(f,g)=(\Phi_y\big|_{\partial D},
\partial\Phi_y\big|_{\partial D})$. 

First, let $v^\varepsilon_y\in\HH(D)={\mathcal S}_D\omega$ 
be an approximation of $v_y$ such that 
\[
\|v_y-v^\varepsilon_y\|_{H^1(D)}<\delta,
\]
where $\delta=\delta(\varepsilon)$ is fixed later. We extend 
$v_y^\varepsilon$
to the whole of $B$ as $v={\mathcal S}_B\omega$.

Having $v^\varepsilon_y$ in $B$, we define $u^\varepsilon_y$ in
$B$ as 
\[
 u^\varepsilon_y = \chi_D u_y + (1-\chi_D)(\Phi_y +v^\varepsilon_y).
\]
We observe that on $\partial B$, we have 
$(\partial u^\varepsilon_y -\partial v
^\varepsilon_y)\big|_{\partial B}=\partial \Phi_y\big|_{\partial B}=0$.
Let us denote 
\[
\partial u^\varepsilon_y\big|_{\partial B} =
\partial v^\varepsilon_y\big|_{\partial B}=L\omega^\varepsilon_y.
\] 
Further,
let us denote $w^\varepsilon_y = T L
\omega^\varepsilon_y$. We show that when 
$\delta$ is small, $w^\varepsilon_y$ and $u^\varepsilon_y$ are close to
each other. To this end, let us define
\[
 r^\varepsilon_y =w^\varepsilon_y-u^\varepsilon_y.
\]
This residual satisfies the equations
\begin{eqnarray*}
\Delta r^\varepsilon_y &=& 0 \mbox{ in $B\setminus\overline D$},\\
\noalign{\vskip4pt}
\nabla\cdot\gamma\nabla r^\varepsilon_y &=& 0 \mbox{ in $D$},
\end{eqnarray*}
the boundary condition
\[
 \partial r^\varepsilon_y\big|_{\partial B} =0,
\]
as well as the transmission conditions
\begin{eqnarray*}
 r^\varepsilon_y\big|_{\partial D}^+ - 
r^\varepsilon_y\big|_{\partial D}^- &=&(v_y-v^\varepsilon_y)
\big|_{\partial D},\\
\noalign{\vskip4pt} 
\partial r^\varepsilon_y\big|_{\partial D}^+ 
-\partial_\gamma r^\varepsilon_y\big|_{\partial D}^+ 
&=& 
(\partial v_y 
-\partial v^\varepsilon_y)\big|_{\partial D}^.
\end{eqnarray*}
By using Green's formula and the trace theorem, 
it is not hard to see that the function
$r^\varepsilon$ satisfies the estimate
\begin{eqnarray*}
 \|r^\varepsilon_y\|_{H^1(B)}&\leq& C(\|v_y-v^\varepsilon_y\|_
 {H^{1/2}(\partial D)}
 +\| \partial v_y 
-\partial v_y^\varepsilon\|_{H^{-1/2}(\partial D)}\\
\noalign{\vskip4pt}
&\leq& C\|v_y-v_y^\varepsilon\|_{H^1(D)}\leq C\delta,
\end{eqnarray*}
and in particular,
\[
 \|r_y^\varepsilon\|_{H^{1/2}(\partial B)}\leq   
C\|r_y^\varepsilon\|_{H^{1}(B)}\leq C\delta.
\]  
Now we have the estimate
\begin{eqnarray*}
\|\Phi_y-(\Lambda-\Lambda_0)L
\omega^\varepsilon_y\|_{H^{1/2}(\partial B)}
&\leq& 
\|\Phi_y-(w^\varepsilon_y -v^\varepsilon_y)\|_{H^{1/2}(\partial B)}\\
\noalign{\vskip4pt}
&\leq& 
\| r^\varepsilon_y\|_{H^{1/2}(\partial B)}\leq C\delta,
\end{eqnarray*}
so by choosing $\delta =\varepsilon/C$ the claim follows.

To prove the second claim of the theorem,
assume that $y\in D$, and let $(v_y,w_y)\in (H^1(D)
\times H^1(D))/\C$ be the solution of the interior transmission
problem with the transmission data $(f,g)=(\Phi_y\big|_{\partial D},
\partial\Phi_y\big|_{\partial D})$.
We show first that as $y$ approaches the boundary $\partial D$, then
$\|v_y\|_{H^1(D)}\rightarrow\infty$. to show this, assume first the
contrary, $\sup_{y\in D}\|v_y\|_{H^1(D)}<\infty$. In particular,
it follows that $\|v_y\|_{H^{1/2}(\partial D)}\leq C$ and 
$\|\partial v_y\|_{H^{-1/2}(\partial D)}\leq C$. We define in $B$
the function $W_y$ as
\[
 W_y=\chi_D w_y +(1-\chi_D)\Phi_y,
\]
satisfying the equations
\begin{eqnarray*}
 \nabla\cdot\gamma\nabla W_y &=& 0 \mbox{ in $D$,}\\
\noalign{\vskip4pt}
 \Delta W_y &=& 0 \mbox{ in $B\setminus\overline D$}
\end{eqnarray*}
with the transmission data
\begin{eqnarray*}
 W_y\big|_{\partial D}^+ - W_y\big|_{\partial D}^- &=&
 -v\big|_{\partial D}^-,\\
\noalign{\vskip4pt}
  \partial W_y\big|_{\partial D}^+ - \partial_\gamma
 W_y\big|_{\partial D}^- &=&
 -\partial v\big|_{\partial D}^-,
\end{eqnarray*}
and the boundary condition
\[
 \partial W_y\big|_{\partial B} =0.
\]
An application of Green's formula leads now to the conlusion
that
\[
 \|W_y\|_{H^1(B)}\leq C(\|v_y\|_{H^{1/2}(\partial D)} + 
\|\partial v_y\|_{H^{-1/2}(\partial D)})\leq C
\]
for all $y\in D$.
In particular, we see that
\[
\sup_{y\in D}\|\Phi_y\|_{H^1(B\setminus D)}<\infty,
\]
which is a contradiction. Thus, we must have 
$\|v_y\|_{H^1(D)}\rightarrow\infty$ as claimed. 

Let $v_y^\varepsilon\in\HH(D)$ be an approximation of
$v_y$ in $H^1(D)$. It follows now that also 
$\|v^\varepsilon_y\|_{H^1(D)}\rightarrow\infty$. By the construction,
$\|L\omega_y^\varepsilon\|_{H^{-1/2}(\partial B)}
\rightarrow\infty$ whis is possible only if
$\|\omega_y^\varepsilon\|_{L^2(\partial \Omega)}
\rightarrow\infty$
as $y$ appraches the boundary $\partial D$. The proof is
complete.
\hfill$\Box$

Finally, let us briefly discuss the case when 
we try to find the
an approximate solution when $y\notin D$. As it is
customary in the
linear sampling approach, we consider the Tikhonov regularized
approximation of the solution to the equation $(\Lambda-\Lambda_0)
\psi_y=\phi_y$. To this end, we need the following result.

\begin{lemma}
Under the assumptions about $\gamma$ made in Section 2,
the operator $\Lambda-\Lambda_0:H^{-1/2}_0(\partial B)
\to H^{1/2}(\partial B)$ is injective and has a dense range.
\end{lemma}

{\em Proof:} To show the injectivity, assume that
$(\Lambda-\Lambda_0)\psi=0$. Set, as usual,
$u=T\psi$ and $v=T_0\psi$, yielding that $w=u-v$
is harmonic in $B\setminus\overline D$ and has
vanishing Cauchy data on $\partial B$. Hence, 
$u=v$ in $B\setminus\overline D$. It follows then
that the pair $(u\big|_D,v\big|_D)$ is a solution
of the interior transmission problem \ref{itp} with vanishing
boundary data, and so Lemma \ref{uniqueness} implies
that $u=v=0$ in $D$ and consequently in the whole of $B$.
Hence, we deduce that also $\psi=0$.

To prove the density, assume the contrary. Then there is
an element $0\neq \eta\in H^{1/2}_0(\partial B)^*$ with
\[
 \langle (\Lambda-\Lambda_0)\psi,\eta\rangle =0
\]
for all $\psi\in H^{-1/2}_0(\partial B)$. Since
$\langle (\Lambda-\Lambda_0)\psi,c\rangle=0$ for
all constants $c$, we may assume that $\eta\in H^{-1/2}_0
(\partial B)$. Further, since $\Lambda-\Lambda_0$ is
symmetric, we deduce that $(\Lambda-\Lambda_0)\eta=0$,
and the injectivity implies $\eta=0$. This contradction
proves the claim.\hfill $\Box$

The above lemma guarantees that we may apply the standard
theory of minimum-norm solutions. In particular (see e.g.
\cite{colton98}), for every $\delta>0$ there is a unique
$\psi^\delta_y\in H^{-1/2}(\partial B)$ that minimizes
the functional
\begin{equation}\label{F}
 F_\alpha(\psi)=\|(\Lambda-\Lambda_0)\psi-\phi_y\|_{H^{1/2}(\partial B)}^2
+\alpha\|\psi\|_{H^{-1/2}(\partial B)}^2,
\end{equation}
with the Morozov discrepancy constraint
\begin{equation}\label{morozov}
\|(\Lambda-\Lambda_0)\psi-\phi_y\|_{H^{1/2}(\partial B)}\leq\delta,
\end{equation}
used to fix the parameter $\alpha=\alpha(\delta)$. By Lemma
\ref{y not in D}, we observe that as $\delta\rightarrow 0+$,
we must have $\|\psi\|_{H^{-1/2}(\partial B)}\rightarrow\infty$. In terms of the regularization
parameter, we have $\alpha\rightarrow 0$ as
$\delta\rightarrow 0$.
What is more, by Lemma \ref{density of L}, for every $\varepsilon>0$
we can always find
an $\omega^{\delta,\varepsilon}_y\in L^2(\partial\Omega)$ such that
\[
 \|(\Lambda-\Lambda_0)(\psi_y^\delta - L\omega^{\delta,\varepsilon}_y)\|
_{H^{1/2}(\partial B)}<\varepsilon.
\]
We can summarize these results in the following theorem that is the
counterpart of Theorem \ref{y in D theorem} when $y\notin D$.

\begin{theorem}\label{y not in D theorem}
Assume that $y\in B\setminus\overline D$ Then for every $\delta>0$
and $\varepsilon>0$ there is an 
$\omega_y^{\delta,\varepsilon} \in L^2(\partial\Omega)$ such that 
\[
\|(\Lambda-\Lambda_0) L\omega^{\delta,\varepsilon}_y-\phi_y\|
_{H^{1/2}(\partial B)}<\delta +\varepsilon,
\]
for which $\|\omega_y^{\delta,\varepsilon}\|_{L^2(\partial\Omega)}
\rightarrow\infty$ as $\delta\rightarrow 0+$.
\end{theorem}

By comparing Theorems \ref{y in D theorem} and \ref{y not in D theorem}
it is not obvious how the linear sampling algorithm should be implemented.
In the articles \cite{colton97a} and \cite{colton00a} (see also the review article \cite{colton00} for further references)
the linear sampling method in inverse scattering
has been studied numerically. Based on those works,
one can suggest the following procedure. 
Given a
'noise level' $\delta>0$, one seeks to minimize 
the functional (\ref{F}) under the constraint
(\ref{morozov}), with
the parameter $y$ varying in a given grid inside
$B$. The norm of the solution
$\psi$ or alternatively, the size of the
regularization parameter $\alpha =\alpha(\delta)$
is used then as a cut-off indicator.

\end{document}